\documentclass{amsart}
\usepackage{amsmath}
\usepackage{framed,comment,enumerate}
\usepackage[pdftex]{graphicx}
\usepackage{epstopdf}
\usepackage{xcolor, framed}
\usepackage{color}

\usepackage{amssymb}
\usepackage{bm}
\usepackage{amsthm}
\usepackage{amscd}
\usepackage{amsmath}
\usepackage{amsfonts}
\usepackage{amssymb}
\usepackage{graphicx}
\usepackage{color}
\newtheorem{theorem}{Theorem}

\newtheorem{lemma}[theorem]{Lemma}

\newtheorem{definition}[theorem]{Definition}

\usepackage{mathrsfs}
\usepackage{titletoc}
\numberwithin{theorem}{section}  
\numberwithin{equation}{section}
\usepackage{arydshln}

\newcommand{\be}{\begin{equation}}
	
	\newcommand{\ee}{\end{equation}}
\newcommand{\bea}{\begin{eqnarray}}
	\newcommand{\eea}{\end{eqnarray}}
\newcommand{\bna}{\begin{eqnarray*}}
	\newcommand{\ena}{\end{eqnarray*}}

\newcommand{\ri}{\right}

\def\R {\mathbb{R}}

\def\PosS{\{u>0\}}
\def\ConS{\{u=0\}}

\def\hem{\hspace{0.5em}}
\def\vem{\vspace{0.6em}}

\theoremstyle{definition}

\numberwithin{equation}{section}

\theoremstyle{remark}

\title[Tangential touch in the Alt-Phillips problem]{Tangential touch between free and fixed boundaries for the fully nonlinear Alt-Phillips problem} 

\author{Yamin Wang}
\address{Department of Mathematics,	National University of Singapore, Singapore}
\email{ym.wang@nus.edu.sg}

\author{Hui Yu}
\address{Department of Mathematics,	National University of Singapore, Singapore}
\email{huiyu@nus.edu.sg}

\begin{document}

\begin{abstract}
For the fully nonlinear Alt-Phillips problem with parameter $\gamma\in(1,2)$, we show that the free boundary intersects the fixed boundary tangentially  where the Dirichlet data vanish. 

For this range of $\gamma$, this result is new even when the operator is the Laplacian.

\end{abstract}
\maketitle

\section{Introduction}
In this article, we address the boundary behavior of a solution to the \textit{fully nonlinear Alt-Phillips problem}. To be precise, for a parameter $\gamma\in(0,2)$ and with the decomposition of coordinates  of $\R^n$ as
$$
x=(x',x_n),
$$ 
 we study solutions to the following system  in the upper half of the unit ball 
\begin{equation}
\label{ap1}
\begin{cases}
F(D^2u)=\gamma u^{\gamma-1}\chi_{\{u>0\}} &\text{ in }B_1\cap\{x_n>0\},\\
u\ge0 &\text{ in }B_1\cap\{x_n>0\},\\
u=0 &\text{ on }B_1\cap\{x_n=0\}.
\end{cases}
\end{equation} 
Here $F$ is a fully nonlinear operator satisfying the following conditions:
\begin{equation} \label{addf1} \begin{cases}
&F \text{ is uniformly elliptic, convex, and  }F(0)=0, \\
&F\in C^1 \text{ with }|DF(M)-DF(0)|\le\omega(|M|) \text{ for a modulus of continuity }\omega.
\end{cases}\end{equation} 

The positive set and the contact set are defined as $\PosS\cap\{x_n>0\}$ and $\ConS\cap\{x_n>0\},$ respectively. Separating them is the \textit{free boundary}, namely, 
\begin{equation} \label{fbd}
\Gamma(u):=\partial\PosS\cap\{x_n>0\}.
\end{equation} 
We are particularly interested in how $\Gamma(u)$ intersects the \textit{fixed boundary} $\{x_n=0\}. $

\vem

The linear counterpart of \eqref{ap1}, that is, when the operator is the Laplacian, is the \textit{classical Alt-Phillips equation}, arising as the Euler-Lagrange equation of the Alt-Phillips functional.  In a domain $\Omega$, this functional is given by
\begin{equation}
\label{EqnAPEnergy}
\mathcal{E}(u;\Omega)=\int_{\Omega}|\nabla u|^2/2+u^\gamma\chi_{\{u>0\}},
\end{equation}
first studied by Phillips \cite{P} and Alt-Phillips \cite{AP}. Apart from its applications in the study of permeable catalysts \cite{A} and population dynamics \cite{GM}, this functional embeds two of the most well-studied free boundary problems as special cases, namely, the Alt-Caffarelli problem when $\gamma=0$ \cite{V} and the obstacle problem when $\gamma=1$ \cite{PSU}. This functional also leads to interesting free boundary problems when the exponent $\gamma$ is in $(-2, 0)$.  As $\gamma\to-2$, one recovers the perimeter functional in the theory of minimal surfaces \cite{DS3, DS4}. 

Motivated by these connections,  the classical Alt-Phillips problem has received intense attention in the past few decades, especially on the interior regularity properties of the solution and of the free boundary. For classical results, the reader may consult \cite{AC, AP, C, CJK, D, DJ, H, W,P}. For modern developments, the reader may consult \cite{DJS,DS1, EdSV,EFeY,FeG, FeY,FiS, JS,KSP, ReRo, SY1, SY2}.

Much less is known about the regularity of the free boundary near the boundary of the domain, even when the operator is the Laplacian. In this direction, the results are known for the cases $\gamma=1$ and $\gamma=0$. 

For $\gamma=1$ and when the operator is the Laplacian,   in a series of works \cite{AMM, AU1, MN,  SU,U},  the authors addressed the boundary behavior of solutions. In particular, they showed that the free boundary intersects the fixed boundary in a tangential fashion where the boundary data vanish. Without this assumption, the situation becomes much more involved \cite{AM, AJ, AS, AU2}.

Still for the problem with the Laplacian operator, when $\gamma=0$, tangential contact between the free boundary and the fixed boundary was first observed by Chang Lara-Savin in \cite{CS}. This leads to the recent breakthrough in Ferreri-Spolaor-Velichkov \cite{FSV}, concerning the boundary branching set. 

\vem

When the operator is fully nonlinear as in our problem \eqref{ap1}, many fundamental tools become inapplicable, most notably, various monotonicity formulae.  This imposes serious challenges on the study of both the interior and boundary regularity of the free boundary. In the interior of the domain, the regularity of the free boundary was addressed for the case $\gamma=0$ in \cite{DFS1, DFS2, W1, W2}, for the case $\gamma=1$ in \cite{L,SY3,B1}, and for $\gamma\in(1,2)$ in \cite{WY}. 

In terms of the behavior near the fixed boundary, so far the only information concerns the case when $\gamma=1$.  For this parameter, Indrei-Minne \cite{IM} and Indrei \cite{I}  showed that the free boundary leaves the fixed boundary  tangentially, following earlier results  in \cite{MM} under more specific assumptions on the operator and the solution. 

\vem
In this work, we study properties of the free boundary near the fixed boundary for the fully nonlinear problem \eqref{ap1}. For $\gamma\in(1,2)$, we establish the tangential contact between the free and the fixed boundaries. 
For this range of $\gamma$, this is new even when the operator is the Laplacian. 

To be precise, our main result reads

\begin{theorem}\label{mt} 
	For $F$ satisfying \eqref{addf1} and $\gamma\in (1, 2)$,  let  $u$ be a solution to \eqref{ap1} with
	$$
0\in\overline{\Gamma (u)}.
	$$
	
	There exists  a  universal modulus of continuity $\sigma$ such that
	\begin{equation*}
		\partial \{u>0\}\cap\{x_n>0\}\subset \{x: x_n\leq \sigma(|x|)|x|\}.
	\end{equation*}
\end{theorem} 
\vem
A constant is called universal if it only depends on the dimension $n$, the exponent $\gamma$, modulus of continuity $\omega$  and the elliptic constants $\Lambda$  (see  \eqref{f1}).

For $\gamma\in (1,2)$, the right-hand side $\gamma u^{\gamma-1}$ degenerates near the free boundary.  Moreover, global solutions are convex for this range of $\gamma$ \cite{WY}. Neither of these holds for $\gamma\in(0,1)$. As result,   we postpone to addressing  the case $\gamma\in (0,1)$ to another work.

\vem

Recall that even when the operator is Laplacian, the boundary behavior of solutions to \eqref{ap1} is only addressed for the case when $\gamma=0$ or $1$.

Compared with the case $\gamma=1$,  when the right-hand side of the equation is constant in the positive set $\PosS$, our equation \eqref{ap1} is more involved.  In particular, when $\gamma=1$, derivatives of the solution solve simple equations in the positive set. This is no longer the case for $\gamma\in  (1,2)$.  

Compared with case $\gamma=0$,  however, our nonlinearity $u^{\gamma-1}$ actually simplifies the situation. For $\gamma=0$, any given constant $a>0$ provides a homogeneous solution to \eqref{ap1} of the form
$$
f_a(x)=a\max\{x_n,0\}.
$$
Moreover, for $a\ge \sqrt{2}$, this is a minimizer to the energy \eqref{EqnAPEnergy}. This multitudes of homogeneous solutions complicated the analysis in \cite{CS, FSV}.

For our problem \eqref{ap1}, however, there is a unique nontrivial solution in the upper half-space. This classification is the heart of our argument (see Section 4).

\vem

This paper is structured as follows. In the next section we provide some preliminary material and introduce some notations.  In Section 3, we discuss the growth of $u$ near the free boundary  $\Gamma(u)$ and its non-degeneracy.  In Section 4,  we classify the solutions to the classical  problem in upper half-space. With these  in Section 5, we complete the proof of Theorem \ref{mt}.

\section{Preliminary}

This section contains three aspects. The first part concerns fully nonlinear elliptic operators. For this,  we refer to \cite{CC}  for a comprehensive treatise.  Secondly  we state some definitions and present a few basic results  related to the fully nonlinear Alt-Phillips problem. The last part is dedicated to the classical Alt-Phillips equation. We collect some preliminary facts about solutions to this equation. 

\subsection{Fully nonlinear elliptic operators.}
Let $\mathcal{S}_n$ denote the space of $n\times n$ symmetric matrices.  Let $\Lambda$ be a constant in $[1, +\infty)$.  A function 
$$
	 F: \mathcal{S}_n \rightarrow \mathbb{R}$$
is a uniformly elliptic operator with ellipticity constant $\Lambda$ if it satisfies 
\begin{equation}\label{f1}
	\frac{1}{\Lambda} \|P\|\leq F(M+P)-F(M)\leq \Lambda \|P\|
\end{equation}
for all $M, P\in \mathcal{S}_n$  and $P\geq 0$. 
 
On top of its ellipticity \eqref{f1}, our assumptions on the operator are
\begin{equation}\label{fy2} F\;\text{is\,convex};\quad  F(0)=0; \end{equation}
 $F$ is $C^1$ with  
\begin{equation}\label{fy1}
|DF(M)-DF(0)|\leq \omega(|M|)
\end{equation}
for a modulus of continuity $\omega$, and 
\begin{equation}\label{fy3}
	DF(0)(M)=trace(M) \hem \text{for\,all}\hem M\in \mathcal{S}_n.	
\end{equation}

A natural concept of weak solution to these operators is that of a viscosity
solution (standard references on the general theory of viscosity solutions include \cite{CC}, \cite{CIL}).
Such solutions enjoy the following stability property. 
\begin{lemma}\label{df} (Stability property)
	Let $(F_j)_{j\geq 1}$  be a sequence of uniformly elliptic operators satisfying \eqref{f1}-\eqref{fy3},
	and $(A_j)_{j\geq 1}$ be a sequence of positive numbers. Let $\Omega$ be an open set in $\mathbb{R}^n$, and  $(u_j)_{j\geq 1}\in C(\Omega)$ be viscosity solutions to
	\begin{equation*}\label{assa1}
		\frac{F_j\left(A_j D^2 u_j\right)}{A_j}=f_j\quad \quad \text{in}\quad \quad \Omega
	\end{equation*}
	for a sequence of continuous functions  $f_j$.  Suppose that as $j\rightarrow \infty$,  	$A_j\rightarrow 0$,
	\begin{equation*}\label{fss1}
		u_j\rightarrow u_{\infty}\quad \quad \text{locally\; uniformly\; in}\quad \Omega,
	\end{equation*}
	and 
	\begin{equation*}\label{ass1}
		f_j\rightarrow f\quad \quad \text{locally\; uniformly\; in}\quad \Omega.
	\end{equation*}

	Then we have 
	$$\Delta u_\infty=f\quad \quad \text{in}\quad \Omega.$$
\end{lemma}	

\proof This follows from classical argument  \cite[Proposition 2.9]{CC} and use of \eqref{fy1} and \eqref{fy3}.  $\hfill\Box$

\subsection{The fully nonlinear Alt-Phillips problem} Recall that we are considering \eqref{ap1} with $F$ satisfying  assumptions \eqref{f1}-\eqref{fy3} above and $\gamma\in (1,2)$.

Let $u$ be a continuous viscosity solution to the above problem. By its local boundedness and  Evans–Krylov estimate   \cite{CC}, we obtain that $u$ satisfies the equation in the classical sense. This implies an overdetermined condition on the free boundary, namely, 
 $$u=|\nabla u|=0\quad \text{on}\quad \Gamma (u).$$

For notational simplicity, we introduce the following notations and classes of solutions.  
\begin{definition}\label{ldv4}  
	Given an open  set $E\subset \mathbb{R}^n$, denote
	$$E^+(\text{or}\hem E_+):= E\cap \{x_n>0\}$$
	and 
	$$E^0(\text{or}\hem E_{0}):= E\cap \{x_n=0\}.$$
\end{definition}

\begin{definition}\label{ld4}  
	For $F$ satisfying \eqref{f1}-\eqref{fy3}, we write that
	$$u\in \mathcal{S}_R^{F}$$
	if $u$ satisfies \eqref{ap1} in $B_R^+$. We write 
	$$u\in \mathcal{S}_{\infty}^{F}$$
	if $u$ satisfies \eqref{ap1} in $\mathbb{R}^n_+$, and such solutions are called global solutions. 
	\vem
	
	An important subclass of $\mathcal{S}_R^{F}$ is 
	$$\mathcal{P}_{R}^{F}(0):=\{u\in \mathcal{S}_R^{F}:\, 0\in\overline{\Gamma (u)} \}$$ 
	where the free boundary $\Gamma(u)$  is defined as in \eqref{fbd}. 
\end{definition}

	With a slight abuse of notation, we use  $\mathcal{S}_R^{\Delta}$ and $\mathcal{S}_{\infty}^{\Delta}$  to denote the classes  when the operator is the Laplacian. \vem

Regarding  the classes of solutions above,  by the maximum principle for the fully nonlinear equations, we  have
\begin{lemma}\label{lf2} (Comparision principle) 
	Let  $u\in \mathcal{S}_1^{F}$.  If $v$ satisfies
\begin{equation*}\label{b6}\begin{cases}
		F (D^2 v) \leq \gamma v^{\gamma -1}&\text{in}\;\;B_1^+, \\
		v\geq 0 &\text{in}\;\; \overline{B^+_1},\\
		v\geq u  &\text{on}\;\; \partial B_1^+,
\end{cases}\end{equation*}
	then 
	\begin{equation*}\label{ass2}
		v\geq u\quad \quad \text{in}\quad  B_1^+.
	\end{equation*}
\end{lemma}

Next we study the symmetry property for the solution to \eqref{ap1}. Let  $u\in \mathcal{S}_R^{F}$ and $x_0\in \overline{\Gamma(u)}$ be a given point. For  $r>0$, we define the rescaled function  on the domain $ D:=\{x\in \mathbb{R}^n: \, x_0+rx\in \overline{B_R^+}\}$,  
\begin{equation}\label{as1}
	u_{x_0, r}(x):=\frac{u(x_0+rx)}{r^{\beta}},
\end{equation}
where 
\begin{equation}\label{asf1}
	\beta:=\frac{2}{2-\gamma}\in (2, +\infty).
\end{equation}
In case $x_0=0$, we usually write for simplicity that 
$$u_{r}:=u_{0, r}.$$
Then   $u_{x_0, r}$ solves \eqref{ap1} in $D$ 
with the nonlinear operator given by 
$$F_r(M):=\frac{1}{r^{\beta-2}}F(r^{\beta-2}M).$$
Notice that $F_r(M)$  still  satisfies the assumptions  \eqref{f1}-\eqref{fy3}.  \vem

In this and next sections,  we frequently  adopt the following terminologies.   Let $(r_j)_{j\geq 1}$ be a sequence of positive numbers such that $r_j\rightarrow 0^+$.  Replacing $r$ by $r_j$ of \eqref{as1},   we consider the sequence of rescalings $u_{x_0, r_j}$,  which is called \textit{blow-up sequence} at $x_0$. The term "\textit{blow-up} of $u$ at $x_0$" will be used for uniform limits of the form  $\lim_{j\rightarrow \infty} u_{x_0, r_j}.$
\vspace{0.04cm}

\subsection{Tools for the  classical  problem}
We gather some results on the Alt-Phillips equation with Laplacian  operator, which are crucial to our study. We will  utilize  the definitions and  terminologies in  Section 2.2.

We begin with the Weiss monotonicity formula that  is essentially due to Weiss \cite{W1}. 
\begin{lemma}\label{l5} (Weiss monotonicity formula) Let $u\in \mathcal{S}_{R}^{\Delta}$, then the function
	\begin{equation*}\label{b1}
		W(u_{ r}, 1)=\int_{B_1^+}\left(\frac{|\nabla u_{r}|^2}{2}+u_{r}^{\gamma}\chi_{\{u_{ r}>0\}}\right)dx-\frac{\beta}{2}\int_{\partial B_1^+\cap \mathbb{R}^n_+}u_{ r}^2d\sigma.
	\end{equation*}
	is non-decreasing in $r\in (0,R)$.  Moreover, $	W(u_{ r}, 1)$ is constant if and only if $u$ is homogeneous of degree $\beta$.  
\end{lemma}

For $\beta$-homogeneous functions, we then present  the  invariance properties of their blow-up limits.  We refer this to  \cite[Lemma 10.9]{V}. 
\begin{lemma}\label{la13} 
Let $u: \, \overline{\mathbb{R}^n_+}\rightarrow \mathbb{R}$ be a $\beta$-homogeneous locally $C^1$ function. Assume  $0\not=p\in \{x_n=0\}$ and $u(p)=0$. 
	Let $r_j\rightarrow 0$ and 
	\begin{equation*}\label{baas32}
		u_{p,\, r_j}(x)=\frac{u(p+r_jx)}{r_j^\beta} 
	\end{equation*} 
	be a blow-up sequence converging locally uniformly to a function $u_{p,\, 0}$ in $ \overline{\mathbb{R}^n_+}$. 
	Then $u_{p,\,0}$ is invariant in the direction $p$, that is 
	$$u_{p, \,0}(x+t p)=u_{p,\, 0}(x) \hem \text{for every } x\in \overline{\mathbb{R}^n_+}\hem \text{and every }  t\in  \mathbb{R}.$$
\end{lemma}

Another useful tool in studying the local structure of the free boundary is the improvement of monotonicity. We prove a version for our problem.  

Roughly speaking, if the solution is `almost' monotone in $B_1^+$ and strictly monotone away from the free boundary, then the results here imply that the solution is indeed monotone  in $B_{1/2}^+$.

  In what follows,  $\{e_1, \dots, e_n\}$ is denoted as the standard basis in $\mathbb{R}^n$ and  $\mathbb{S}^{n-1}$ is the $(n-1)$ sphere $\{x\in \mathbb{R}^n: |x|=1\}$.

\begin{lemma}\label{l12} (Improvement of monotonicity)
	Let $u\in \mathcal{S}_{1}^{\Delta}$.  Suppose that for some direction  $e\in \mathbb{S}^{n-1}_+$ and for some positive constants $\epsilon, \eta, \delta$, we have 
	$$\partial_e u\geq -\epsilon\quad \quad \text{in} \quad B_1^+$$ 
	and 
	\begin{equation}\label{b17}
		\partial_e u\geq  \delta>0\quad \quad \text{in}\quad B_1^+\cap \{u\geq \eta >0\}.
	\end{equation}

Then 	there exists universal constant $c_0>0$ such that if 
	$$\eta\leq c_0\quad \text{and} \quad  \epsilon\leq c_0 \delta,$$ 
	we have $$\partial_e u\geq 0 \quad \quad \text{in} \quad B_{1/2}^+.$$
\end{lemma}
\proof In view of  \eqref{b17} and $\partial_e u=0$ on $\{u=0\}\cap B_{1/2}^+$,  it suffices to prove $\partial_e u\geq 0$ for  $\{0<u<\eta\}\cap B_{1/2}^+$. To do so, we will use a barrier. Set the domain 
$$\Omega=\{0<u<\eta\}\cap B^+_{3/4}.$$
Pick  $x_0\in B^+_{1/2}\cap \{0<u<\eta\}$, and define the  barrier function on $\overline{\Omega}$ as 
$$h:=\partial_e u+\delta (|x-x_0|^{\beta}-\widetilde{C}_0u),$$
where the constant  $\widetilde{C}_0>>\beta(n+\beta-2)$ that will be determined later. For $u\in  \mathcal{S}_{1}^{\Delta}$,   define the linearized operator as
\begin{equation}\label{bq16}
	Lw:=\Delta w-\gamma(\gamma-1)u^{\gamma-2}w\quad \text{in}\quad \Omega,
\end{equation}
then one has
\begin{equation}\label{b16}
	Lh=\delta\left(\beta(n+\beta-2)|x-x_0|^{\beta-2}-\gamma(\gamma-1)u^{\gamma-2}|x-x_0|^{\beta}-\widetilde{C}_0\gamma(2-\gamma)u^{\gamma-1}\right). 
\end{equation}

First we estimate the value of $h$ on the boundary of $\Omega$. On $\{u=\eta\}\cap B^+_{3/4}$ and $\{u=\eta\}\cap \partial B^+_{3/4}$, 
\begin{equation*}
	h\geq \delta (1-\widetilde{C}_0\eta)>0
\end{equation*}
due to \eqref{b17} and
\begin{equation}\label{be16}
	\eta\leq c_0:= \frac{1}{\widetilde{C}_0\,2^{2(\beta+1)}}<\frac{1}{\widetilde{C}_0}<\frac{1}{\beta(n+\beta-2)}.
\end{equation}
Along $\{u=0\}\cap B^+_{3/4}$, from  $\partial_e u=0$, one finds
\begin{equation*}\label{b19}
	h\geq \delta |x-x_0|^{\beta}\geq 0.
\end{equation*}
On $\{u=0\}\cap \partial B^+_{3/4}$,  applying $\langle e, e_n \rangle \geq 0$ yields  $\partial_e u\geq 0$ and thus
\begin{equation*}\label{b20}
	h\geq \partial_e u \geq 0.
\end{equation*}
On $ \partial B^+_{3/4} \cap \{0<u<\eta\}$, since $\epsilon\leq c_0\delta$ and $\eta\, \widetilde{C}_0  \leq 1/2^{2(\beta+1)}$, with \eqref{be16}, 
\begin{equation*}\label{b21}
	h\geq -\epsilon+\delta \left(\frac{1}{2^{2\beta}}- \widetilde{C}_0  \eta\right)\geq\left(  \frac{3}{2^{2(\beta+1)}} -\frac{1}{\widetilde{C}_0\,2^{2(\beta+1)}}\right) \delta>0.
\end{equation*}
To summarize, 
\begin{equation}\label{b18}
	h\geq 0\quad \quad \text{along}\quad \partial \Omega.
\end{equation} 

Next we prove that $Lh\leq 0$ in $\Omega$. In fact, from \eqref{b16}, if 
$$\beta(n+\beta-2)|x-x_0|^{\beta-2}\leq   \gamma(\gamma-1)u^{\gamma-2}|x-x_0|^{\beta},$$ 
then the desired result follows immediately. For otherwise, one has
$$C_1 u^{2-\gamma}\geq |x-x_0|^2,\quad \quad\quad C_1:= \frac{\beta(n+\beta-2)}{\gamma(\gamma-1)},$$
this means that
$$ |x-x_0|^{\beta-2}\leq (C_1\,u)^{\gamma-1}.$$
Since $\widetilde{C}_0>>\beta(n+\beta-2)$, we choose 
$$\widetilde{C}_0 \geq  \frac{\beta(n+\beta-2) C_1^{\gamma-1}}{\gamma(2-\gamma)}$$ to get 
$$\beta(n+\beta-2)|x-x_0|^{\beta-2}\leq \widetilde{C}_0\gamma(2-\gamma)u^{\gamma-1}.$$ 
Consequently, 
\begin{equation}\label{b22}
	Lh\leq 0 \quad \quad \text{in} \quad \Omega.
\end{equation}
Combining \eqref{b18} and \eqref{b22}, by the comparison principle, we obtain
$$h\geq 0\quad \text{in}\quad \{0<u<\eta\}\cap B_{3/4}^+,$$
and thus 
$$ \partial_e u \geq \delta \widetilde{C}_0 u(x_0)\geq 0 \quad \quad \text{in} \quad  \{0<u<\eta\}\cap B_{1/2}^+.$$ 
This finishes the proof of the lemma. $\hfill\Box$\vem 

If $u\in \mathcal{S}_{1}^{\Delta}$ is monotone in the direction $e$,   our  next lemma shows that either $u$ is strictly monotone  in its  positive set   or $u$ is a trivial solution.  Precisely, 

\begin{lemma}\label{lq4} 
	Let $u\in \mathcal{S}_{1}^{\Delta}$.  If for some $e\in \mathbb{S}^{n-1}_+$, 
	$$\partial_e u\geq 0\quad \text{in}\quad B_1^+,$$
	then either
	$$
	\partial_e u>0\quad \text{in}\quad \{u>0\}\cap B_1^+,
	$$
	or 
	$$
	u\equiv0\quad \text{in}\quad B_1^+.
	$$
\end{lemma}
\proof Suppose not.  There exists a point $p\in \{u>0\}\cap B_1^+$ such that 
\begin{equation}\label{cas1}	\partial_e u (p)=0,\end{equation}
but with
\begin{equation}\label{cas3}
	u\not\equiv0\quad \quad \text{in}\quad B_1^+.
\end{equation}

For $u\in \mathcal{S}_{1}^{\Delta}$, we define the linearized operator as in \eqref{bq16}.
Since 
\begin{equation*} \begin{cases}
		L(\partial_e u)=0&\text{in}\;\{u>0\}\cap B_1^+,\vspace{0.1cm} \\
		\partial_e u \geq 0&\text{in}\;\;B_1^+,
\end{cases}\end{equation*}
we derive from   \eqref{cas1} and the  maximum principle that
$$	\partial_e u \equiv 0\quad \quad \text{in}\quad \{u>0\}\cap B_1^+.$$
Therefore, 
$$	\partial_e u \equiv 0\quad\quad  \text{in}\quad B_1^+.$$
This together with $u=0$ on $B_1^0$  gives 
$$u(x+te)-u(x)=\int_{0}^t 	\partial_e u dt=0$$
for $t>0$ and $x\in B_1^0$.  As a consequence,  
$$	u\equiv0\quad\quad  \text{in}\quad B_1^+,$$
a contradiction to \eqref{cas3}.
$\hfill\Box$\vem

Before ending this section, we show one more lemma, as observed by Alt-Phillips \cite{AP} and  Bonorino \cite{B2}.  

\begin{lemma}\label{l4} 
	Let $u\in \mathcal{S}_{R}^{\Delta}$. Set $w=u^{\frac{2}{\beta}}$. Then $w$ solves
	$$\Delta w=A-B\frac{|\nabla w|^2}{w}\quad \quad \text{on}\quad \{w>0\}\cap B_R^+,$$ 	
	where $A=\gamma(2-\gamma)$ and $B=\frac{\gamma-1}{2-\gamma}$. Moreover, the solution $w$ is $C^{1, 1}_{loc}( B_R^+)$. \vem
\end{lemma}

\section{The growth and nondegeneracy of solutions } 
In this section, assuming $0\in \overline{\Gamma(u)}$,  we mainly discuss the growth estimates and  nondegeneracy of the solution  $u\in \mathcal{P}_R^{F}(0)$. Here recall that the free boundary   $\Gamma(u)$  and the set $\mathcal{P}_R^{F}(0)$ are respectively denoted by  \eqref{fbd} and Definition \ref{ld4}.  \vem

We start with an observation on the  gradient of $u$ at the points where  the free boundary  leaves the fixed boundary.

\begin{lemma}\label{l6} Assume $u\in \mathcal{P}_1^{F}(0)$, then 
	$$u(0)=|\nabla u(0)|=0.$$
\end{lemma}
\proof It suffices to prove $\partial_{x_n} u(0)=0$. Argue  by contradiction that $\partial_{x_n} u(0)>0$. By $C^{1, \alpha}$  continuity at $x=0$, we have $\partial_{x_n} u(0)>0$ in $B_r^+$ for small $r>0$, contradicting $0\in\overline{\Gamma(u)}$. Thus  $\partial_{x_n} u(0)=0$ and $|\nabla u(0)|=0$.  $\hfill\Box$\vem

Then we  establish $\beta$-power growth control of solutions near the free boundary.  In the following, 
$B_r(x_0)\subset \mathbb{R}^n$ denotes  the open ball of radius $r$ centered at $x_0$, and recall $B_r^+(x_0)=B_r(x_0)\cap \{x_n>0\}$.
\begin{lemma}\label{l7} (The growth of solutions) Let $u\in \mathcal{P}_1^{F}(0)$.  There is a universal constant $C>0$ such that
	$$\sup_{B_r^+(x_0)} u\leq Cr^{\beta}$$ 
	for every  $x_0\in \overline{\Gamma(u)\cap B^+_{1/2}}$ satisfying   
	\begin{equation}\label{caf1}
		u(x_0)=|\nabla u(x_0)|=0
	\end{equation}
	and for  every $r\in (0, 1/2)$. 
\end{lemma}
\proof  It is enough to prove that there exists a universal constant $C>0$ such that for every  $x_0\in \overline{\Gamma(u)\cap B^+_{1/2}}$ with \eqref{caf1} and $j\in \mathbb{N}$,
\begin{equation}\label{ca1}
	\sup_{B^+_{2^{-j}(x_0)}}u\leq C(2^{-j})^{\beta}.
\end{equation}
For this purpose, we define the set
$$\mathbb{M}(u, x)=\Big\{j\in\mathbb{N}\;| \sup_{B^+_{2^{-(j+1)}}(x_0)} u \geq  2^{-\beta} \sup_{B^+_{2^{-j}}(x_0)}u \Big\}.$$ 
Observe that if $\mathbb{M}$ is empty then one may  obtain the desired estimate by iteration.

We suppose $\mathbb{M} \not =\emptyset$.  If \eqref{ca1} fails, then there exists a sequence of solutions  $u_j\in \mathcal{P}_1^{F}(0)$ with $x_j\in \overline{\Gamma(u_j)\cap B^+_{1/2}}$  satisfying  \eqref{caf1},  and $k_j\in \mathbb{N}$ ($k_j\rightarrow \infty$) such that
\begin{equation}\label{c1}
	\sup_{B^+_{2^{-(k_j+1)}}(x_j)} u_j \geq  2^{-\beta} \sup_{B^+_{2^{-k_j}}(x_j)}u_j,\quad  \quad 
	\sup_{B^+_{2^{-k_j}}(x_j)}u_j >j (2^{-k_j})^{\beta}.
\end{equation}
Now denote $ \tilde{u}_j$ as
\begin{equation}\label{cq1}
	\tilde{u}_j(y)=\frac{u_j(x_j+2^{-k_j}y)}{\sup_{B^+_{2^{-(k_j+1)}}(x_j)} u_j }\quad \text{in}\quad B_1(0)\cap \{y_n\geq -l_{j} \},
\end{equation}
where for simplicity we write 
$$l_{j}:= \frac{(x_j)_n}{2^{-k_j}}.$$
From above, it is obvious to see 
\begin{equation}\label{b9}
	\sup \tilde{u}_j(y)=1 \quad \text{in}\quad  B_{1/2}(0)\cap \{y_n\geq -l_{j}\}
\end{equation}
and according to  \eqref{c1}, 
\begin{equation}\label{b10}
	\tilde{u}_j(y)\leq 2^\beta \quad \text{in}\quad  B_{1}(0)\cap \{y_n\geq -l_{j}\}.
\end{equation}

Next we introduce
\begin{equation}\label{bs8}
	\widetilde{F}_j(M)=\frac{2^{-2k_j}}{\sup_{B^+_{2^{-(k_j+1)}}(x_j)} u_j }
	F_j\left(\frac{\sup_{B^+_{2^{-(k_j+1)}}(x_j)} u_j}{2^{-2k_j}}M\right)
\end{equation}
Due to \eqref{c1}, \eqref{cq1} and   \eqref{b10}, in $B_1(0)\cap \{y_n\geq -l_{j}\}$, there holds 
\bea\nonumber
0 \leq 	\widetilde{F}_j(D^2  \tilde{u}_j(y))&=&\frac{2^{-2k_j}}{\sup_{B^+_{2^{-(k_j+1)}}(x_j)} u_j }F_j(D^2 u_j (x_j+2^{-k_j}x))\\\nonumber
&\leq&\frac{(2^{\beta})^{2-\gamma}}{j^{2-\gamma} }\gamma (\tilde{u}_j (y))^{\gamma-1}\\\label{b8}
&\leq &\frac{\gamma 2^{\beta}}{j^{2-\gamma}}\rightarrow 0\quad \text{as}\quad j\rightarrow \infty.	\eea

Combining  with \eqref{b10}, \eqref{b8} and $\tilde{u}_j=0$ on $\{y_n=-l_{j}\}$, by the boundary $C^{1, \alpha}$-regularity (see e.g. \cite{LZ, SS}), we obtain
$$
\| \tilde{u}_j \|_{C^{1, \alpha}(B_{1/2}(0)\cap \{y_n\geq -l_n\})}\leq C,
$$
where $0<\alpha<1$ and $C$ is a universal constant. Thus, up to a subsequence (still denoted by $\tilde{u}_j$),  
\begin{equation}\label{bcd10}
	\tilde{u}_j\rightarrow \tilde{u}_{\infty}\quad \quad \text{in}\quad C^{1, \, \alpha}_{loc}(B_{1/2}(0)\cap \{y_n\geq -l_n\})
\end{equation}
as $ j\rightarrow \infty$.\vem

To derive the equation $\tilde{u}_{\infty}$ solves,  from  the boundary $C^{2, \alpha}$ regularity \cite{LZ, SS},  we know that
\begin{equation}\label{bc10}
	\| u_j \|_{C^{2, \alpha}(\overline{B_{3/4}^+(0)})}\leq C.
\end{equation}
This  with  $x_j\in \overline{\Gamma(u_j)\cap B^+_{1/2}}$ implies
 $$Du_j\left(x_j\right)=0,\quad \quad D^2u_j\left(x_j\right)=0,$$ 
 where we used Lemma \ref{l6} if $x_j=0$. 
Due to the above and  \eqref{bc10}, 
\begin{equation*}\label{bacx10} \sup_{B^+_{2^{-(k_j+1)}}(x_j)} u_j \leq C 2^{-(k_j+1)(2+\alpha)}.\end{equation*}
It then follows that 
\begin{equation}\label{bcw10}
	\frac{\sup_{B^+_{2^{-(k_j+1)}}(x_j)} u_j}{2^{-2k_j}}\rightarrow 0\quad \text{as}\quad j\rightarrow \infty.
\end{equation}
By \eqref{bs8}-\eqref{bcd10} and \eqref{bcw10}, using the stability properties (see Lemma \ref{df}),  we conclude  that $\tilde{u}_{\infty}$ solves
\be\label{b11}\left\{\begin{array}{lll}
	\tilde{u}_{\infty}(y)\geq 0	\quad\;\; \;\text{in}\;\;B_{1/2}(0)\cap \{y_n\geq -\lim_{j\rightarrow \infty }\sup l_{j}\},\vspace{0.15cm} \\
	\Delta \tilde{u}_{\infty}(y) = 0\quad \text{in}\;\;B_{1/2}(0)\cap \{y_n\geq -\lim_{j\rightarrow \infty }\sup l_{j}\},\vspace{0.15cm} \\
	\tilde{u}_{\infty}(0)=|\nabla \tilde{u}_{\infty}(0)|=0.
\end{array}\ri.\ee
From here we need to separate two cases depending on the limit of $l_{j}$ as  $j\rightarrow \infty$. 	\vem

Case 1. If $l_{j}\rightarrow 0$, note that  $\tilde{u}_{\infty}\geq 0$ and $\tilde{u}_{\infty}\not\equiv 0$ on $B_{1/2}(0)\cap\{y_n\geq 0\}$.  Since $\tilde{u}_{\infty}(0)=0$,  by Hopf's Lemma, we reach a contradiction with $|\nabla \tilde{u}_{\infty}(0)|=0$.\vem

Case 2. If $l_{j}\geq \epsilon$ for some small $\epsilon>0$, applying   the maximum principle to (\ref{b11}), one has $\tilde{u}_{\infty}(x)\equiv 0$ in $B_{1/2}(0)\cap \{y_n\geq -\lim_{j\rightarrow \infty }\sup l_{j}\}$. However this contradicts \eqref{b9}. \vem

Hence the lemma is established.
$\hfill\Box$\vem

With the help of Lemmas \ref{l6}, \ref{l7} and regularity theory \cite{LZ, SS},  the following   lemma provides us with enough compactness for a family of rescaled solutions.

\begin{lemma}\label{l8} ($C^{1, \alpha}$ regularity up to the boundary) Let $u\in\mathcal{P}_{1}^{F}(0)$ and  $u_r$ be  defined in \eqref{as1}.
	One has $$\|u_r\|_{C^{1, \alpha}(\overline{B_{1/2}^+})}\leq C\quad \quad \text{for}\quad \alpha\in (0,1).$$ 
\end{lemma}

This lemma is important if we want to blow up (see the definition in Section 2.2) the solution  $u\in \mathcal{P}_{1}^{F}(0)$ near the free boundary.  It assures that the blow up sequence is uniformly bounded on compact sets, and thus the limit of blow-up is bounded on compact sets.  

However, it is  not clear whether  $u_r$ is not the zero function as $r\rightarrow 0$. If $u(x)=o(|x|^{\beta})$,  then any $u_r$  as $r\rightarrow 0$ is identically zero. To prevent this, we need a nondegeneracy from below. 

\begin{lemma}\label{l10} (Nondegeneracy)
	Let $u\in\mathcal{P}_{1}^{F}(0)$. 	If  $x_0\in \overline{\{u>0\}}$,  then there is a universal constant $C>0$ such that 
	\begin{equation}\label{bc9}
		\sup_{\partial B_r^+(x_0)} u \geq Cr^{\beta}
	\end{equation}
	for every $r>0$ such that  $B_r^+(x_0)\subset B_1^+(0)$.
\end{lemma}
\proof Let us define 
$$  \Phi (x) := a |x|^{\beta} \quad  \;\, \text{in}\quad \overline{B_1^+(0)},$$
where the constant $a$ is small such that $a^{\gamma-2}\geq \frac{\Lambda \beta^2}{2}$. 	\vem

If \eqref{bc9} is false,  there exists  $x_0\in \overline{\{u>0\}}$  and a ball $B_r^+(x_0)\subset B_1^+(0)$ such that
$$u<c_0 \Phi (x-x_0)\quad \text{on} \quad \partial B_r^+(x_0)$$ 
for positive constant $c_0$. On one hand, from above,  for $p\in B_r^+(x_0)$ with $|p-x_0|<<1$, 
\begin{equation}\label{bch8}
	u(x)\leq \Phi(x-p)\quad \text{on}\quad \partial B_r^+(x_0).
\end{equation}
On the other hand, 
\begin{equation}\label{bc8}
	\begin{split}
		F(D^2 \Phi)\leq \Lambda a \beta(\beta-1) |x|^{\beta(\gamma-1)}\leq \gamma \Phi^{\gamma-1}\quad \text{in}\quad B_1^+(0).
	\end{split}
\end{equation}

Exploiting \eqref{bch8}, \eqref{bc8} and  Lemma \ref{lf2}, we know
$$u(x)\leq \Phi (x-p)\quad \text{in}\;\;B_r^+(x_0).$$ 
Therefore $u(p)\leq \Phi(0)=0$, contradicting $x_0\in \overline{\{u>0\}}$. $\hfill\Box$\vem

In conclusion,   assume   $u\in \mathcal{P}_{1}^{F}(0)$, then by Lemma \ref{l8} and Arzal$\grave{\mathrm{a}}$-Ascoli theorem, there exists a sequence $r_j\rightarrow 0$ such that 
\begin{equation}\label{bce8}
	u_{r_j}\rightarrow u_0\quad \quad \text{in}\quad C^{1, \alpha}_{loc}(\overline{\mathbb{R}^n_+}),
\end{equation}
where $u_0\not\equiv0$ due to Lemma \ref{l10} and $u_0\in \mathcal{S}_{\infty}^{\Delta}$ by Lemma \ref{df}.  In particular,  $u_0$ also satisfies $\beta$-growth and nondegeneracy properties as in Lemmas 
\ref{l7} and \ref{l10}.\vem

\section{Classification of global  solutions to the classical problem}
In this section, we will use the technical tools developed in Section 2.3
to classify the global solution $u_0\in C^{1, \alpha}_{loc}(\overline{\mathbb{R}^n_+})$ of the classical  Alt-Phillips equation. In what follows, we  omit the subscript  of $u_0$ for simplicity. 

Concretely,   given $ u$ satisfying
\be\label{d1}\left\{\begin{array}{lll}
	u \in C^{1, \alpha}_{loc}(\overline{\mathbb{R}^n_+}),\vspace{0.1cm} \\
	\Delta u=\gamma u^{\gamma-1} 	\quad \text{in}\;\;\mathbb{R}^n_+,\vspace{0.1cm} \\
	u=0 \quad\quad\quad\;\;\;\,\text{on}\;\;\mathbb{R}^n_{0}, \vspace{0.1cm} \\
	u\geq 0\quad\quad\quad\;\;\;\; \text{in}\;\;\mathbb{R}^n_+, 
\end{array}\ri.\ee
and  \vspace{0.15cm}\\
\indent  (i)\; $aR^\beta\leq \sup_{\partial B_R^+(x_0)}u$ for $ x_0\in \overline{\{u>0\}},\;\forall R>0$,\vspace{0.15cm}\\
\indent   (ii)\; $AR^\beta \geq \sup_{B_R^+(x_0)}u $ for  $ x_0\in \overline{\Gamma(u)} $ with $ u(x_0)=|\nabla u(x_0)|=0,\; \forall R>0$,\vspace{0.15cm} \\
where $a$ and $A$ are positive universal constants, the main result of this section  reads 

\begin{lemma}\label{class} 
	Let $u$ be a solution to  \eqref{d1}. Assume the assumptions (i)-(ii) hold, then 
	\begin{equation}\label{cs1}
		u(x)=\left(\sqrt{2}\left(\frac{x_n}{\beta}\right)^+\right)^{\beta}.
	\end{equation}
\end{lemma}

To understand the configuration of the  solution to \eqref{d1}, 
we shall reduce the problem to the study of homogeneous
solutions by a blow-down procedure. The proof of Lemma \ref{class} is postponed to  the end of this section.

\begin{lemma}\label{l11} ($\beta$-homogeneity of solutions)
	Let $u$ be a solution to  \eqref{d1}. Assume the assumptions (i)-(ii) hold.  For a sequence $R_j\rightarrow\infty$, define
	\begin{equation}\label{bcr8}
		{u}_{R_j}(x)=\frac{1}{R_j^{\beta}}u(R_j x). 
	\end{equation}
	Then, passing to a subsequence (still denoted by ${u}_{R_j}$), we have
	$${u}_{R_j}\rightarrow u_{\infty} \quad \; \text{in} \quad C^{1, \alpha}_{loc}(\overline{\mathbb{R}^n_+}),$$ 
	where $u_{\infty}$ is a nonzero $\beta$-homogeneous solution to  \eqref{d1}.  
\end{lemma}
\proof  From  the assumption (ii) and regularity theory, we know that ${u}_{R_j}\in  C^{1, \alpha}_{loc}(\overline{\mathbb{R}^n_+})$.  The convergence to  solution  $u_{\infty}$  of  \eqref{d1}  follows directly from  $C^{1, \alpha}$ regularity for ${u}_{R_j}$ and the nondegeneracy of $u_{\infty}$   is  obtained by  assumption (i). 

Now we prove the homogeneity of $u_{\infty}$. For any $\rho>0$, we infer from Lemma \ref{l5}  that 
$$W(u_{\infty}, \rho)=\lim_{R_j\rightarrow \infty}W({u}_{R_j}, \rho)=\lim_{R_j\rightarrow \infty}W(u,\, \rho R_j)=W(u, \infty).$$
Here the last equality exists due to the monotonicity of $W$ and since $u\in C^{1, \alpha}_{loc}(\overline{\mathbb{R}^n_+})$,  this gives $W(u, R)\leq C$ for all $R>0$ and for some  constant $C>0$. Thus  $u_{\infty}$ is $\beta$ homogeneous. $\hfill\Box$\vem

On order to  proceed, we require a technical  lemma  as follows. 
\begin{lemma}\label{l13} 
	Assume that $u$ is a $\beta$-homogeneous solution to  \eqref{d1}.  Under assumption (i), there holds
	\begin{equation}\label{bd31}
		\partial_{x_n} u\geq 0\quad \text{in}\quad \overline{\mathbb{R}^n_+}.
	\end{equation}
\end{lemma}
\proof By induction, suppose that \eqref{bd31} holds for dimensions $1, \cdots, n-1$. Now we prove (\ref{bd31}) is true in $n$-dimension.  In view of  the homogeneity of $u$, it suffices to prove 
\begin{equation}\label{b40}
	\partial_{x_n} u\geq 0   \quad \text{in} \quad  \partial B_1^+. 		
\end{equation}
The proof is  split into  the following  two steps. 
\vem

\emph{Step 1.} In this step, we shall prove $\partial_{x_n}u\geq 0$ near the fixed boundary $ \partial B_1^+\cap \{x_n= 0\}$.   We claim that there is a universal constant $r_0>0$ such that 
\begin{equation}\label{b31}
	\partial_{x_n} u(x)\geq 0\quad\text{in}\quad \partial B_1^+\cap \{0\leq x_n\leq r_0\}. \end{equation}
To this end, we need to decompose  three cases. \vem

Case 1. If  there is a point $p\in \partial B_1^+\cap \{ x_n=0\}$ satisfying $\partial_{x_n} u(p)>0$, since $u\in C^{1, \alpha}_{loc}(\overline{\mathbb{R}^n_+})$, we deduce
$$
\partial_{x_n} u(p)\geq 0\quad \text{in}\quad \partial B_1^+\cap \{0\leq x_n\leq r_1(p)\},
$$
where the constant  $r_1(p)>0$  depends on $p$. Then by the compactness, we can find  a universal $r_0$ such that the claim (\ref{b31}) holds.  \vem

Case 2. If $p\in \partial B_1^+\cap \{ x_n=0\}$ with $p\not\in \{\overline{u>0}\}$, then $\partial_{x_n} u(p)=0$. By the $C^{1, \alpha}$ regularity of $u$, there exists $r_1(p)$ such that
$$
\partial_{x_n} u(p)= 0\quad \text{in}\quad \partial B_1^+\cap \{0\leq x_n\leq r_1(p)\}.
$$
Then (\ref{b31}) follows with no difficulty. \vem

Case 3. For $p\in \partial B_1^+\cap \{ x_n=0\}$ with $p\in \{\overline{u>0}\}$ and $\partial_{x_n} u(p)=0$, we set
\begin{equation*}\label{ba32}
	u_{p,\, r}(x)=\frac{u(p+rx)}{r^\beta}. 
\end{equation*} 
Assume  with no loss of generality that $p=e_1\in  \overline{\Gamma(u)}$. 
By the $C^{1, \alpha}$ regularity and $\beta$-homogeneity of   $u$,    
\begin{equation}\label{bw32}
	u_{e_1, r_j}\rightarrow u_{e_1, 0} \quad \text{in}\quad C^{1, \alpha}_{loc}( \overline{\mathbb{R}^n_+})
\end{equation} 
along a subsequence $r_j\rightarrow 0$.  Since  $u(e_1)=0$, employing   Lemma  \ref{la13}  gives
\begin{equation*}\label{b32}
	u_{e_1, 0}(x+t e_1)=u_{e_1, 0}(x) \quad \forall x\in \overline{\mathbb{R}^n_+},\;\; t\in  \mathbb{R}.
\end{equation*} 
This means that in $\overline{\mathbb{R}^n_+}$,
$$\partial_{e_1} u_{e_1, 0}\equiv 0.$$  From  the initial hypothesis,  as a result, 
\begin{equation}\label{bq32}
\partial_{x_n} u_{p, 0}\geq 0\quad \text{in}\quad \overline{\mathbb{R}^n_+}. \end{equation}

We assert that there exist positive constants $\delta$ and $\eta\leq c_0$  ($c_0$ is given by  Lemma \ref{l12}) such that
\begin{equation}\label{ba37}
\partial_{x_n} u_{p, r_j}>\delta\quad \quad \text{in}\quad \{u_{p, r_j}>\eta\} \cap B_1^+(p).
\end{equation}
Indeed, with Lemma \ref{lq4}  and \eqref{bq32}, we have 
\begin{equation}\label{baa37}
	\partial_{x_n} u_{p, 0}>0\quad  \quad \text{in}\quad \{u_{p, 0}>0\}\cap B_1^+(p).
\end{equation}
Otherwise, $u_{p,\, 0}\equiv 0$ in $B_1^+$. But  in light of the assumptions (i), 
$$\sup_{\partial B_1^+(p)} u_{p, r_j} \geq a>0,$$
which leads to a contradiction. Then by \eqref{baa37},  there is $\eta\leq c_0$ and  $\delta>0$ such that
\begin{equation*}
	\partial_{x_n} u_{p, 0}>2\delta\quad \text{in}\quad \left\{u_{p, 0}\geq \frac{\eta}{2} \right\}\cap B_1^+(p).
\end{equation*}
This combining wth the convergence \eqref{bw32} yields  \eqref{ba37}, as asserted.

Meanwhile, in view of  \eqref{bw32} and  \eqref{bq32}, one can find a  positive constant $\epsilon\leq c_0\delta$ to get
\begin{equation}\label{ba36}
	\partial_{x_n} u_{p, r_j}\geq -\epsilon\quad \text{in} \quad B_1^+(p)
\end{equation}
for $j$ sufficiently large. 

Lastly together with \eqref{ba37} and  \eqref{ba36}, invoking   Lemma \ref{l12} concludes 
$$	\partial_{x_n} u_{p, r_j}\geq 0\quad \quad \text{in}\quad B^+_{1/2}(p).$$
Equivalently $\partial_{x_n} u\geq 0$ in $B^+_{r_j/2}(p)$. Hence \eqref{b31} is claimed. \vem

\emph{Step 2.} In this step,  we  shall complete \eqref{b40}. To this end,  let 
$$v=u^{2/\beta}.$$ 
With Step 1,  it is enough to verify $v_n\geq 0$ in $\partial B_1^+\cap \{x_n>r_0\}$, where 
$$v_n:=\partial_{x_n} v. $$  
Utilizing  Lemma \ref{l4}, one easily get 
\begin{equation}\label{b37}
	\Delta v_n=B \,v_n \frac{|\nabla v|^2}{v^2}  -2B\frac{\nabla v\cdot \nabla v_n }{v}\quad \text{in}\quad \{v>0\}\cap \partial B_1^+
\end{equation}
with
$$B:=\frac{\gamma-1}{2-\gamma}>0.$$

We argue by contradiction that there exists a positive constant $m$ such that 
\begin{equation}\label{b38}
	\inf_{\partial B_1^+\cap \{v>0\}}\frac{v_n }{x_n}=-m<0. 
\end{equation}
Let $r_0$ be as in \eqref{b31}.  For $x\in \partial B_1^+\cap \{v>0\}$,  if $0<x_n\leq r_0$,  by Step 1, there holds $v_n\geq 0$. Thus from \eqref{b38}, we have 
$$x_n> r_0.$$

By the definition of infimum, there exists a subsequence $p_j\in \partial B_1^+\cap \{v>0\}$ with $(p_j)_n> r_0$ such that
$$\frac{v_n(p_j) }{(p_j)_n}\rightarrow -m. $$
In view of  the compactness, 
$$p_j\rightarrow p_\infty\in \partial B_1^+\cap \overline{ \{v>0\}}$$
along with
\begin{equation}\label{d38}
	\frac{v_n(p_\infty) }{(p_\infty)_n}  =-m.   
\end{equation}
We have two possibilities to consider. 

If $p_{\infty}\in \{v>0\}\cap\partial B_1^+$,  we  define  a function $\varphi$ on $ \{v>0\}\cap \partial B_1^+$ as
$$\varphi=v_n+mx_n.$$
Note that $\varphi$ is one homogenous.  From \eqref{b38} and \eqref{d38}, one can see 
\begin{equation}\label{b39}
	\varphi\geq 0\quad \text{in}\quad B_1^+, \quad \quad \varphi(p_{\infty})=0.
\end{equation}
It then follows that  $\nabla v_n=-me_n$ at $p_{\infty}$. This together with  \eqref{b37} yields
\begin{equation*}
	\begin{split}
		\Delta \varphi  (p_{\infty})&=-B\frac{mx_n }{v^2}|\nabla v|^2 +2B\frac{\nabla v\cdot me_n }{v}\\
		&= -B\frac{mx_n }{v^2}|\nabla v|^2 -2B\frac{m^2x_n}{v}<0, 
	\end{split}
\end{equation*}
contradicting  \eqref{b39}.  

If $p_{\infty}\in \partial\{v>0\}\cap \partial B_1^+$, since $(p_{\infty})_n\geq r_0$,  one has $v_n(p_{\infty})=0$. However, 
$$v_n(p_{\infty})=-m(p_{\infty})_n<-m r_0,$$
contradiction again.  \vem 

This confirms \eqref{b40} and completes the proof of the lemma.    $\hfill\Box$\vem

Relying on the above analysis, the next lemma gives a  classification result for $\beta$-homogeneous solution. 
\begin{lemma}\label{l14} 
	Let $u$ be a $\beta$-homogeneous solution to (\ref{d1}). Under assumption (i), we have  \begin{equation}\label{fek}
		u(x)=u(x_n)=\left(\sqrt{2}\left(\frac{x_n}{\beta}\right)^+\right)^{\beta}.
	\end{equation}
\end{lemma}
\proof Recall $\mathbb{S}^{n-1}_+=\mathbb{S}^{n-1}\cap\{x_n> 0\}$. Define 
$$\mathcal{C}:=\left\{e\in \mathbb{S}^{n-1}_+\,|\, \partial_eu(x)\geq 0\; \text{for\,all} \; x\in \overline{\mathbb{R}_+^n}\right\}.$$
Notice  that $\mathcal{C}$  is closed and by  Lemma \ref{l13},  $\mathcal{C}\not=\emptyset$. Assume $e^{\ast}\in \partial \mathcal{C}$ with $(e^{\ast})_n>0$. Then    there exists  a sequence $e_k\in \mathcal{C}$ converging to $e^{\ast}$ with $ \partial_{e_k}u\geq 0$ as $k\rightarrow \infty$. Therefore, $ \partial_{e^{\ast}}u\geq 0$ and $e^{\ast} \in \mathcal{C}$.

For a small constant $t>0$,  let 
$$e_t:=\frac{e^{\ast}+te_1}{|e^{\ast}+te_1|}.$$
Obviously $\langle e_t, e_n \rangle > 0$. We  estimate  
\bea\nonumber
\partial_{e_t} u=\nabla u\cdot e_t&=&
\frac{1}{|e^{\ast}+te_1|}\left( \partial_{e^{\ast}}u+t \partial_{e_1}u\right)\\\label{b50} 
&\geq& \frac{1}{1+t}\partial_{e^{\ast}} u-\frac{t}{1-t}|\partial_{e_1} u|. 
\eea

Since $e^{\ast} \in \mathcal{C}$, by Lemma \ref{lq4} and assumption (i),  we derive 
\begin{equation*}\label{bs50} \partial_{e^{\ast}}  u>0\quad \text{in}\quad \{u>0\}\cap B_1^+.
\end{equation*} 
Thus there exists $\eta\leq c_0$ (with $c_0$ in  Lemma \ref{l12}) and $\delta>0$ such that
\begin{equation}\label{b53} 
	\partial_{e^{\ast}} u\geq 2 \delta  \quad \quad \text{in}\quad \{u\geq \eta\}\cap B_1^+.
\end{equation}
Note that $|\partial_{e_1} u|\leq C$ in $B_1^+$. Putting this and  \eqref{b53} into \eqref{b50}, we obtain
\begin{equation}\label{b52} 
\partial_{e_t} u\geq \delta  \quad \quad \text{in}\quad \{u\geq \eta\}\cap B_1^+
\end{equation}
provided that $t$ is small enough.

On the other hand,  since $\partial_{e^{\ast}}u\geq 0$,  by  \eqref{b50},  one can easily check 
\begin{equation}\label{b51} 
\partial_{e_t} u\geq -\frac{t}{1-t}|\partial_{e_1} u |\geq -\epsilon \; \quad \text{in}\quad  B_1^+
\end{equation}
where $\epsilon\leq c_0\delta$.

Thanks to Lemma \ref{l12}, we conclude by \eqref{b52} and \eqref{b51} that  
$$\partial_{e_t} u\geq 0\quad \quad \text{in}\quad B_{1/2}^+.$$
Thus $e_t\in \mathcal{C}$, contradicting $e^{\ast}\in \partial  \mathcal{C}$. Consequently  $\partial_{e_i}u\equiv 0$ for $i=1, 2, \dots, n-1$. So 
$u(x)=u(x_n)$  and the problem \eqref{d1} is reduced to solving an ODE.  Then \eqref{fek}  follows immediately. $\hfill\Box$\vem

\emph{Proof of Lemma \ref{class}.} Let $u$ be the  solution to \eqref{d1} and $u_{R_j}$ be defined by \eqref{bcr8}.  Observe by Lemma \ref{l11} that 
\begin{equation}\label{bvy57}
	{u}_{R_j}\rightarrow u_\infty\hem\hem  \text{in}\hem C^{1, \alpha}_{loc}(\overline{\mathbb{R}^n_+})
\end{equation}
with  $u_\infty$ given by  \eqref{fek}.  Applying  Lemma \ref{l14},  for any direction $e\in \mathbb{S}^{n-1}_+$,   one has
\begin{equation*}\label{bvy56}
	\partial_{e} u_\infty(x)\geq c_1\eta^{1-1/\beta}\langle e, e_n\rangle>0\hem\hem \text{in}\hem \left\{u_\infty\geq\frac{\eta}{2}\right\}\cap B_1^+,
\end{equation*}
where $c_1>0$  is universal constant and $\eta\leq c_0$ with $c_0$ in Lemma \ref{l12}.  With this and  \eqref{bvy57},  we  deduce 
\begin{equation}\label{bv56}
\partial_{e}  u_{R_j}(x)\geq \frac{c_1}{2}\eta^{1-1/\beta}\langle e, e_n\rangle\geq \delta>0\hem \hem \text{in}\hem  \left\{u_{R_j}>\eta\right\}\cap B_1^+.		
\end{equation}
Moreover, since   $\partial_{e} u_\infty(x)\geq 0$ in $B_1^+$,   by  \eqref{bvy57},  there exists $\epsilon\leq c_0\delta$ such that 
\begin{equation}\label{cx1}
\partial_{e} u_{R_j}\geq -\epsilon\quad \text{in} \hem B_1^+.
\end{equation}

Now exploiting  \eqref{bv56}, \eqref{cx1} and  Lemma \ref{l12}, we arrive at 
$$\partial_{e} u_{R_j}\geq 0\quad \text{in} \hem B_{1/2}^+.$$
In particular, 
$$\partial_{e} u\geq 0\hem \hem  \text{in}\hem B_{R_j/2}^+,\quad \forall R_j>>0.$$
As a result,
$$\partial_{e} u\geq 0\hem \hem  \text{in}\hem\mathbb{R}^n_+,$$
then repeating the proof as in Lemma \ref{l14} gives  \eqref{cs1}, as desired. $\hfill\Box$

\section{Proof of Theorem \ref{mt}}
\emph{Completion of the proof of Theorem \ref{mt}.} Let $u\in \mathcal{P}_{1}^{F}(0)$ (see Definition \ref{ld4}). It is enough to prove that for any $\epsilon>0$, there is $\rho_\epsilon>0$ such that
$$\Gamma(u)\cap B_{\rho_\epsilon}^+\subset B_{\rho_\epsilon}^+\backslash A_\epsilon,\quad \quad \quad A_\epsilon:=\{x_n>\epsilon|x'|\}.$$ 
By way of contradiction,  there exists $u_j\in \mathcal{P}_{1}^{F}(0)$ and $\epsilon>0$ such that for all $j\in \mathbb{N}$,
$$x_j\in \Gamma(u_j)\cap B^+_{1/j}\cap A_\epsilon.$$
Observe that $r_j=|x_j|\rightarrow 0$ as $j\rightarrow \infty$.  Consider  rescalings
$$\tilde{u}_j(x):=\frac{u_j(r_j x)}{r_j^\beta}, $$
where $\beta$ is in (\ref{asf1}). 
For each function $\tilde{u}_j(x)$, one has a point 
$$\tilde{x}_j=\frac{x_j}{r_j}\in \Gamma(\tilde{u}_j)\cap \partial B_1^+\cap A_\epsilon.$$
Then by compactness (see Lemma \ref{l8}), over a subsequence, 
$$\tilde{u}_j \rightarrow u_0\in \mathcal{S}_{\infty}^{\Delta}\hem \hem  \text{in}\hem C^{1, \alpha}_{loc}(\overline{\mathbb{R}^n_+})$$
with 
\begin{equation*}\label{b55} 
	\tilde{x}_j \rightarrow x_0\in \partial B_1^+\cap \overline{A}_\epsilon.
\end{equation*}
In particular, we have by nondegeneracy 
\begin{equation}\label{b54} 
	x_0\in \Gamma(u_0). 
\end{equation}

However,  by  (\ref{bce8}) and Lemma  \ref{class}, the blow-up of $u$ at the origin is  $$u_0(x)=\left(\sqrt{2}\left(\frac{x_n}{\beta}\right)^+\right)^{\beta}.$$
In particular $\Gamma(u_0)=\emptyset$, contradicting \eqref{b54}.  The proof of the theorem is thus completed. $\hfill\Box$\vem


\end{document}